\documentclass[12pt, leqno]{amsart}
\usepackage[mathscr]{eucal}
\usepackage{amssymb,amsmath}
\usepackage{amsfonts}
\usepackage{a4}

\newtheorem{theorem}{Theorem}

\newtheorem{lemma}{Lemma}
\newtheorem{corollary}{Corollary}

\begin{document}
\title{On splitting polynomials with noncommutative coefficients}

\author[Tomasz Maszczyk]{Tomasz Maszczyk\dag}
\address{Institute of Mathematics\\
Polish Academy of Sciences\\
Sniadeckich 8\newline 00--956 Warszawa, Poland\\
\newline Institute of Mathematics\\
University of Warsaw\\ Banacha 2\newline 02--097 Warszawa, Poland}

\email{t.maszczyk@uw.edu.pl}

\thanks{\dag The author was partially supported by KBN grants N201 1770 33 and 115/E-343/SPB/6.PR UE/DIE 50/2005-2008.}

\begin{abstract} It is shown that for every splitting of a polynomial with noncommutative coefficients into
linear factors $(X-a_{k})$ with $a_{k}$'s commuting with
coefficients, any cyclic permutation of linear factors gives the
same result and all  $a_{k}$ are roots of that polynomial. It
implies that although the set of $a_{k}$'s appearing in a splitting of a
polynomial with commutative coefficients in some noncommutative
extension does not determine the splitting (in general), the cyclic
order consisting of roots appearing in the splitting does. Examples of this phenomenon are given.

\end{abstract}

\maketitle

\paragraph{\textbf{1. Introduction}} Let $f(X)=f_{n}X^{n}+f_{n-1}X^{n-1}+\cdots + f_{0}\in
A[X]$ be a polynomial with coefficients in a commutative unital
ring $A$. Suppose there is given a splitting of $f(X)$ in $A[X]$
\begin{align}
f(X)=f_{n}(X-a_{1})\cdots(X-a_{n}).
\end{align}
Then by the substitution homomorphism argument one sees that all
$a_{k}$'s are roots of $f(X)$ and by commutativity of $A[X]$ any
permutation of them defines the same splitting. Therefore the
problem of splitting of a given polynomial reduces to the problem
of finding the set of its roots. This fact is fundamental for
Galois theory and algebraic geometry.

In the case of noncommutative coefficients of a given polynomial
the situation is much worse. First of all, a given splitting does
not reduces to the set of elements $a_{k}$, since we cannot
permute linear factors because of noncommutativity of $A[X]$.
Moreover, if $a\in A$ is not central in $A$ then the substitution
homomorphism of rings
\begin{align}
\mathbb{Z}[X]\rightarrow A,\ \ X\mapsto a
\end{align}
does not extend to a homomorphism of $A$-algebras
\begin{align}
A[X]\rightarrow A,\ \ X\mapsto a,
\end{align}
because $X$ is central in $A[X]$. This means that one can not use
the substitution $A$-algebra homomorphism argument to prove that
elements $a_{1}, \ldots, a_{n}$ appearing in the decomposition
\begin{align}
f(X)=f_{n}(X-a_{1})\cdots(X-a_{n})
\end{align}
are roots of $f(X)$. The problem of such splittings in terms of
relationships between coefficients of a given polynomial with a
generic set of its (left or right) roots and elements $a_{k}$ (so
called pseudoroots) was related to quadratic algebras with
structure encoded by graphs in
\cite{gr}\cite{grw}\cite{ggrw}\cite{rsw}. However, these
relationships are much more complicated than in the commutative
case and make sense only if some elements of the algebra are
invertible.

The interest for splitting polynomials in noncommutative algebras
started in 1921 when Wedderburn proved \cite{w} that any minimal
polynomial $f(X)\in K[X]$ of an element of a central division
algebra $A$ algebraic over the center $K$ of $A$ splits in $A[X]$
into linear factors which can be permuted cyclically  and every
pseudoroot appearing in this splitting is a root of $f(X)$. This
fact was very helpful in determining the structure of division
algebras of small order \cite{w} and found many other applications
(see e.g. \cite{ll} for references).

The aim of this paper is to show that under the assumption that
coefficients (which do not have to commute one with each other)
commute with pseudoroots (which do not have to commute one with
each other) the situation is much closer to the commutative case.
We show that then pseudoroots are roots and any cyclic permutation
of them gives the same splitting. This means that instead of
finite sets of commutative roots (ordered $n$-tuples up to all
permutations) we obtain finite cyclically ordered sets (ordered
$n$-tuples up to all cyclic permutations) of noncommutative roots.

We give examples of such splittings, where in spite of cyclic
symmetry, transposition of any two consequtive linear factors is
impossible. These examples are minimal in the sense that the
algebra in which we split our polynomial is generated by the roots
appearing in the splitting.

It does not seem that this basic (elementary, trivially provable
etc but not trivial at all) fact could be easily derived from the
known theory of splitting polynomials in noncommutative algebras
(Gelfand-Retakh-Wilson \cite{grw}). The reason is that the
condition of commutativity between coefficients and pseudoroots is
a closed condition while the general theory works for generic
elements. We will see in examples that our result is true even if
differences of pseudoroots are not invertible or even nilpotent.

\vspace{3mm}
\paragraph{\textbf{2. Results}}
\begin{lemma}
Let $g, h$ be elements of a monoid $G$, where $h$ is right
cancelable and the product $gh$ commutes with $h$. Then $g$ and
$h$ commute.
\end{lemma}
\emph{Proof.} Since $gh$ commutes with $h$, we have
\begin{align}
ghh=hgh.
\end{align}
But $h$ is right cancelable, hence consequently
\begin{align}
gh=hg.\ \ \Box
\end{align}

\begin{lemma}
Let $A$ be a  ring and $A^{a}$ be its subring of elements
commuting with a fixed $a\in A$. If $f(X)\in A^{a}[X]$ decomposes
in $A[X]$ as follows
\begin{align}
f(X)=g(X)(X-a)\label{decomp}
\end{align}
then $g(X)\in A^{a}[X]$ and $f(a)=0$.
\end{lemma}

\emph{Proof.} To prove that $g(X)\in A^{a}[X]=A[X]^{X-a}$ take
$G=A[X]$ with multiplication of polynomials, $g=g(X)$,
$h=h(X)=X-a$ and apply Lemma 1.

To prove the root property apply to (\ref{decomp}) the
substitution homomorphism
\begin{align}
A^{a}[X] & \rightarrow A^{a},\\
     X& \mapsto a.\ \ \Box\nonumber
\end{align}

\begin{theorem}
Let $A$ be a unital ring and $A^{a_{1},\ldots, a_{n}}$ be its
subring of elements commuting with $a_{1}, \ldots, a_{n}\in A$. If
$f(X)\in A^{a_{1},\ldots, a_{n}}[X]$ splits in $A[X]$ as follows
\begin{align}
f(X)=f_{n}(X-a_{1})(X-a_{2})\cdots(X-a_{n})
\end{align}
then
\begin{align}
f(X)=f_{n}(X-a_{n})(X-a_{1})\cdots(X-a_{n-1})
\end{align}
and
\begin{align}
f(a_{1})=\cdots = f(a_{n})=0.
\end{align}
\end{theorem}

\emph{Proof.} To prove the cyclic property of the splitting take
$a=a_{n}$, $g(X)=f_{n}(X-a_{1})(X-a_{2})\cdots(X-a_{n-1})$ and
apply Lemma 2. Then Lemma 2 also implies the root property for
$a_{n}$. By the cyclic property the same holds for other
$a_{k}$'s. $\Box$

\begin{corollary}
For any splitting as in Theorem 1 substitution homomorphisms
\begin{align}
A^{a_{1},\ldots, a_{n}}[X] & \rightarrow A^{a_{1},\ldots, a_{n}}[a_{k}]\subset A^{a_{k}}\subset A\\
X & \mapsto a_{k}\nonumber
\end{align}
define a ring homomorphism
\begin{align}
A^{a_{1},\ldots, a_{n}}[X]/(f(X))\rightarrow A^{a_{1},\ldots,
a_{n}}[a_{1}]\times\cdots\times A^{a_{1},\ldots, a_{n}}[a_{n}]
\end{align}
where cyclic permutations of roots in the splitting correspond to
cyclic permutations of factors in the cartesian product on the
right hand side.
\end{corollary}

\paragraph{\textbf{3. Examples}}\ \

\vspace{3mm}\emph{Example 1}. Let $A$ be the ring of upper
triangular $2\times 2$ matrices over a nonzero commutative ring
$K$ and take $f(X)=X^{2}(X-1)\in K[X]\subset A[X]$. Although it
has a double root in $K$ it can be split in $A$ as follows
\begin{align}
f(X)=(X-a_{1})(X-a_{2})(X-a_{3}),
\end{align}
where
\begin{align}
a_{1}=\left(\begin{array}{cc}
             0 & 0\\
             0 & 1
             \end{array}\right),\ \  a_{2}=\left(\begin{array}{cc}
             0 & -1\\
             0 & 0
             \end{array}\right),\ \  a_{3}=\left(\begin{array}{cc}
             1 & 1\\
             0 & 0
             \end{array}\right)
\end{align}
are pairwise distinct roots in $A$. This can be viewed as a kind of resolution of singularity by passing to a noncommutative extension.

The linear factors can be cyclically
permuted, but none two of them can be transposed, because
\begin{align}
[a_{1}, a_{2}]= [a_{2}, a_{3}]=[a_{3}, a_{1}]=-a_{2}\neq 0.
\end{align}

Since the Vandermonde matrix
\begin{align}
\left(\begin{array}{ccc}
             a_{1}^{2} & a_{2}^{2} & a_{3}^{2}\\
             a_{1} & a_{2} & a_{3}\\
             1 & 1 & 1
             \end{array}\right)=\left(\begin{array}{cccccc}
             0 & 0 & 0 & 0 & 1 & 1 \\
             0 & 1 & 0 & 0 & 0 & 0 \\
             0 & 0 & 0 & -1 & 1 & 1 \\
             0 & 1 & 0 & 0 & 0 & 0 \\
             1 & 0 & 1 & 0 & 1 & 0 \\
             0 & 1 & 0 & 1 & 0 & 1
             \end{array}\right)
             \end{align}
is not invertible, we are far from the theory presented in
\cite{ggrw}. However, it is a remarkable fact, that this example
shares many properties with splitting of a polynomial in its Galois
extension.

First of all, it is obvious that $A^{a_{1}, a_{2}, a_{3}}=K$. Moreover, $A$
is freely generated by $a_{1}, a_{2}, a_{3}$ as a module over $K$ with a
multiplication table

\vspace{3mm} \centerline{\begin{tabular}{|l||l|l|l|}\hline
                 & $a_{1}$                   & $a_{2}$       &  $a_{3}$    \\
\hline\hline
$a_{1}$          & $a_{1}$                   & 0             &  0          \\
\hline
$a_{2}$          & $a_{2}$                   & 0             &  0          \\
\hline
$a_{3}$          & $-a_{2}$                  &  $a_{2}$      &  $a_{3}$    \\
\hline
\end{tabular}}

\vspace{3mm}
 If $K$ doesn't contain nontrivial idempotents all
endomorphisms of the extension $K\subset A$ come in families
$\varepsilon, \varepsilon', \varepsilon^{\sigma}_{s},
\varepsilon_{s}$ parameterized by elements $\sigma$ and $s$, where
$\sigma$'s are elements of the multiplicative monoid of $K$ acting
(from the right) on elements $s$ of the (right) $K$-module $K$ by
right multiplication. The logic of this notation will be clear
later, when the rules of matrix multiplication
\begin{align}
\left(\begin{array}{cc}
             \tau & 0\\
             t & 1
             \end{array}\right)
             \left(\begin{array}{cc}
             \sigma & 0\\
             s & 1
             \end{array}\right)=
             \left(\begin{array}{cc}
             \tau\sigma & 0\\
             t\sigma+s & 1
             \end{array}\right),\ \
             (t, 1)\left(\begin{array}{cc}
             \sigma & 0\\
             s & 1
             \end{array}\right)=(t\sigma+s, 1)
             \end{align}
will appear in the structure of the endomorphism monoid and its
various actions.

The endomorphism monoid is determined by its values on basic
elements $a_{1}, a_{2}, a_{3}$ as follows

\vspace{3mm}
\centerline{\begin{tabular}{|l||l|l|l|}\hline
                   & $a_{1}$                 & $a_{2}$      &  $a_{3}$                 \\
\hline\hline
$\varepsilon$           & 1                       & 0            & 0                        \\
\hline
$\varepsilon'$            & 0                       & 0            &  1                       \\
\hline
$\varepsilon^{\sigma}_{s}$ & $a_{1}+sa_{2}$          & $\sigma a_{2}$ &  $(1-\sigma-s)a_{2}+a_{3}$ \\
\hline
$\varepsilon_{s}$       & $(1-s)a_{2}+a_{3}$ & 0            &  $a_{1}+sa_{2}$         \\
\hline
\end{tabular}}
\vspace{3mm}\begin{flushleft}hence the monoid structure is
\end{flushleft}

\vspace{3mm} \begin{center}\begin{tabular}{|l||l|l|l|l|}\hline
                     & $\varepsilon$ & $\varepsilon'$       &        $\varepsilon^{\tau}_{t}$       & $\varepsilon_{t}$      \\
\hline\hline
$\varepsilon$             & $\varepsilon$ & $\varepsilon$      &   $\varepsilon$                     & $\varepsilon$          \\
\hline
$\varepsilon'$              & $\varepsilon'$  & $\varepsilon'$       &   $\varepsilon'$                      & $\varepsilon'$           \\
\hline
$\varepsilon^{\sigma}_{s}$ & $\varepsilon$ & $\varepsilon'$       & $\varepsilon^{\tau\sigma}_{t\sigma+s}$ & $\varepsilon_{t}$      \\
\hline
$\varepsilon_{s}$         & $\varepsilon'$  & $\varepsilon$      & $\varepsilon_{t\sigma+s}$             &  $\varepsilon^{0}_{t}$  \\
\hline
\end{tabular}\end{center}

\vspace{3mm} \begin{flushleft}with the neutral element
$\varepsilon^{1}_{0}$. Note that the only invertible elements are $\varepsilon^{\tau}_{t}$ with $\tau$ invertible in $K$.\end{flushleft}

If $K$ is a domain all roots  and all cycles of $f(X)$ also come
in families parameterized by elements $\tau$ and $t$.

The families of roots are
\begin{eqnarray*}
r^{\tau} & = &\tau a_{2}, \\
r^{\tau}_{t} & = &(1-\tau-t)a_{2}+a_{3}, \\
r_{t}  & = & a_{1}+ta_{2}, \\
r  & = & a_{1}+ a_{2}+a_{3}.
\end{eqnarray*}

The  families of cycles are
\begin{eqnarray*}
c^{\tau} & = & (r^{\tau}, r^{-\tau}, r),\\
c^{\tau}_{t} & = & (r^{\tau}, r^{\tau}_{t}, r_{t}),\\
c_{t}  & = & (r^{0}, r_{t},  r^{0}_{t}).
\end{eqnarray*}

\begin{flushleft}We see that\end{flushleft}

(a) The set of all roots is the union of supports of all cycles.
\begin{flushleft}One can check that\end{flushleft}

(b) The only cycles whose roots form a basis of $A$ as a free $K$-module are $c^{\tau}_{t}$ with $\tau$ invertible in $K$. 

\vspace{3mm}The action of the endomorphism monoid on roots is

\vspace{3mm} \begin{center}\begin{tabular}{|l||l|l|l|l|}\hline
                     & $r^{\tau}$ & $r^{\tau}_{t}$       &        $r_{t}$       & $r$      \\
\hline\hline
$\varepsilon$             & $r^{0}$   & $r^{0}$             &   $r$                     & $r$          \\
\hline
$\varepsilon'$              & $r^{0}$  & $r$       &   $r^{0}$                      & $r$           \\
\hline
$\varepsilon^{\sigma}_{s}$ & $r^{\tau\sigma}$ & $r^{\tau\sigma}_{t\sigma  +s}$       & $r_{t\sigma  +s}$ & $r$      \\
\hline
$\varepsilon_{s}$         & $r^{0}$  & $r_{s}$      & $r^{0}_{s}$             &  $r$  \\
\hline
\end{tabular}\end{center}

\begin{flushleft}which implies that\end{flushleft}

(c) Every root is a translate of $r^{1}$ or $r^{1}_{0}$.

(d) Every root can be translated to $r^{0}$ or $r$.

\vspace{3mm} The action of the endomorphism monoid on cycles is

\vspace{3mm} \centerline{\begin{tabular}{|l||l|l|l|}
\hline
                           & $c^{\tau}$                & $c^{\tau}_{t}$       &  $c_{t}$    \\
\hline\hline
$\varepsilon, \varepsilon'$          & $c^{0}$                   & $c^{0}$              &  $c^{0}$    \\
\hline
$\varepsilon^{\sigma}_{s}$      & $c^{\tau\sigma}$          & $c^{\tau\sigma}_{t\sigma+s}$   &  $c_{t}$          \\
\hline
$\varepsilon_{s}$               & $c^{0}$                  &  $c_{s}$             &  $c^{0}_{s}$    \\
\hline
\end{tabular}}

\vspace{3mm}

\begin{flushleft}which implies that\end{flushleft}

(e) Every cycle is a translate of $c^{1}$ or $c^{1}_{0}$.

(f) Every cycle can be translated to $c^{0}$.

\begin{flushleft}Points (e) and (f) replace the fact that the Galois group of a splitting field of a separable polynomial permutes the roots. Points (c) and (d) replace the fact that this action is transitive.\end{flushleft}

Note however that in our case there are many cycles and endomorphisms can move roots from one cycle to another. It turns out that roots are no more equivalent. Instead of strict equivalence of roots induced by the transitive Galois action we have the action of endomorphisms on some lattice of roots. This can be described as follows.

\vspace{3mm}
The lattice of ideals in $K[X]$
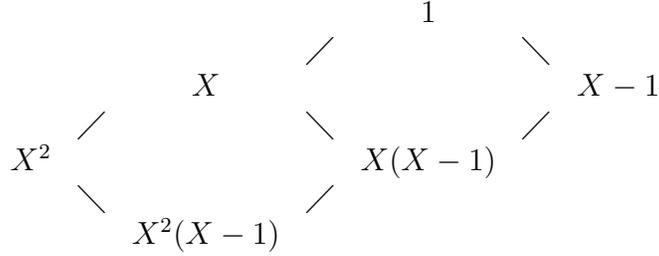
$$\begin{array}{ccccccc}
 & & & & 1 & & \\
 & & & \diagup & & \diagdown & \\
 & & X & & & & X-1 \\
 & \diagup & & \diagdown & & \diagup & \\
 X^{2} & & & & X(X-1) & & \\
 & \diagdown & & \diagup & & & \\
 & & X^{2}(X-1) & & & &
\end{array}$$
defines a partial order on roots by taking the minimal polynomial

\vspace{3mm} \begin{center}\begin{tabular}{|l|l|l|l|l|l|}\hline
      $r^{0}$     & $r^{\tau}$, $\tau\neq 0$ & $r^{\tau}_{t}$       &        $r_{t}$       & $r$      \\
\hline\hline
$X$             & $X^{2}$   & $X(X-1)$             &   $X(X-1)$                     & $X-1$          \\
\hline
\end{tabular}\end{center}

\vspace{3mm}
\begin{flushleft}which looks as follows\end{flushleft}
$$\begin{array}{ccccccc}
 & & r^{0} & & & & r \\
 & \diagup & & \diagdown & & \diagup & \\
 r^{\tau},\ \tau\neq 0 & & & & r^{\tau}_{t}, r_{t} & &
\end{array}$$

One can check that all endomorphisms act along the above poset of roots moving them at most upwards (it is obvious that no endomorphism can diminish a root with
respect to the partial order induced by the minimal polynomial).

Among all endomorphisms only $\varepsilon$ and $\varepsilon'$ do not preserve this
partial order (there are essentially four
exceptions: $r^{0}\succ r^{\tau}_{t}$, $r^{0}\succ r_{t}$, $r\succ
r^{\tau}_{t}$, $r\succ r_{t}$ but
$\varepsilon'(r^{0})\nsucc\varepsilon'(r^{\tau}_{t})$,
$\varepsilon(r^{0})\nsucc\varepsilon(r_{t})$,
$\varepsilon(r)\nsucc\varepsilon(r^{\tau}_{t})$,
$\varepsilon'(r)\nsucc\varepsilon'(r_{t})$), although they do preserve
the weaker partial order opposite to that which is induced by the degree (levels in the poset of roots). In particular, all automorphisms do preserve the above partial order on roots. 
This replaces transitivity of the Galois action on roots.

\vspace{3mm}\emph{Example 2}. Let $A$ be the ring of $3\times 3$
matrices over a nonzero commutative ring $K$ and take
$f(X)=X^{3}-4\in K[X]\subset A[X]$. It can be split as follows
\begin{align}
f(X)=(X-a_{1})(X-a_{2})(X-a_{3}),
\end{align}
with three pairwise distinct roots in $A$ 
\begin{align}
a_{1}=\left(\begin{array}{ccc}
             0 & 2 & 0\\
             0 & 0 & 2\\
             1 & 0 & 0
             \end{array}\right),\ \  a_{2}=\left(\begin{array}{ccc}
            0 & -1 & 0\\
             0 & 0 & 2\\
             -2 & 0 & 0
             \end{array}\right),\ \  a_{3}=\left(\begin{array}{ccc}
             0 & -1 & 0\\
             0 & 0 & -4\\
             1 & 0 & 0
             \end{array}\right).
\end{align}
Note that after reduction modulo 3 $f(X)\equiv (X-1)^{3}$ and $a_{1}\equiv a_{2}\equiv a_{3}$. Therefore this can be viewed as a kind of resolution of singularity in positive characteristic by passing to a noncommutative extension of a lifting to rational integers.
Again, the linear factors can be cyclically permuted. If $3\neq 0$
in $K$ none two of them can be transposed, because
\begin{align}
[a_{1}, a_{2}]= [a_{2}, a_{3}]=[a_{3},
a_{1}]=\left(\begin{array}{ccc}
             0 & 0 & 6\\
             -6 & 0 & 0\\
             0 & 3 & 0
             \end{array}\right)\neq 0.
\end{align}

Moreover, if 2 and 3 are not zero divisors in $K$, again
$A^{a_{1}, a_{2}, a_{3}}=K$. Indeed, $A^{a_{1}, a_{2},
a_{3}}\subset A$ consists of elements
\begin{align}
a=\left(\begin{array}{ccc}
             \alpha & \beta & 2\gamma\\
             -2\gamma & \alpha & \beta\\
             -\beta & \gamma & \alpha
             \end{array}\right),
             \end{align}
where $3\beta=0$, $6\gamma=0$.

If $2,3,5$ and $7$ are invertible in $K$, $A$ is generated by
$a_{1}, a_{2}, a_{3}$ as an algebra over $K=A^{a_{1}, a_{2},
a_{3}}$. Indeed, then we have

\begin{align}
\left(\begin{array}{ccc}
             1 & 0 & 0\\
             0 & 0 & 0\\
             0 & 0 & 0
             \end{array}\right) & =-\frac{1}{9}a_{1}^{2}a_{2}-\frac{1}{18}a_{2}^{2}a_{3},
             \left(\begin{array}{ccc}
             0 & 1 & 0\\
             0 & 0 & 0\\
             0 & 0 & 0
             \end{array}\right)
             =-\frac{2}{3}a_{2}+\frac{1}{15}a_{1}^{2}a_{2}^{2}-\frac{1}{10}a_{2}^{2}a_{3}^{2},\nonumber\\
\left(\begin{array}{ccc}
             0 & 0 & 1\\
             0 & 0 & 0\\
             0 & 0 & 0
             \end{array}\right) & =-\frac{1}{3}a_{1}^{2}+\frac{5}{6}a_{2}^{2}+a_{3}^{2},
             \left(\begin{array}{ccc}
             0 & 0 & 0\\
             1 & 0 & 0\\
             0 & 0 & 0
             \end{array}\right)
             =\frac{1}{2}a_{1}^{2}-a_{2}^{2}-a_{3}^{2},\nonumber\\
\left(\begin{array}{ccc}
             0 & 0 & 0\\
             0 & 1 & 0\\
             0 & 0 & 0
             \end{array}\right) & =\frac{1}{3}a_{1}a_{2}^{2}+\frac{5}{14}a_{1}^{2}a_{2}+\frac{2}{21}a_{2}^{2}a_{3},
             \left(\begin{array}{ccc}
             0 & 0 & 0\\
             0 & 0 & 1\\
             0 & 0 & 0
             \end{array}\right)=\frac{2}{9}a_{1}+\frac{2}{9}a_{2}-\frac{1}{36}a_{1}^{2}a_{2}^{2},\nonumber\\
\left(\begin{array}{ccc}
             0 & 0 & 0\\
             0 & 0 & 0\\
             1 & 0 & 0
             \end{array}\right) & =-\frac{1}{2}a_{2}-\frac{1}{20}a_{1}^{2}a_{2}^{2}-\frac{1}{20}a_{2}^{2}a_{3}^{2},
             \left(\begin{array}{ccc}
             0 & 0 & 0\\
             0 & 0 & 0\\
             0 & 1 & 0
             \end{array}\right)=\frac{2}{3}a_{1}^{2}-\frac{2}{3}a_{2}^{2}-a_{3}^{2},\nonumber\\
\left(\begin{array}{ccc}
             0 & 0 & 0\\
             0 & 0 & 0\\
             0 & 0 & 1
             \end{array}\right) & =\frac{1}{6}a_{1}a_{2}^{2}-\frac{1}{9}a_{1}^{2}a_{2}-\frac{2}{9}a_{2}^{2}a_{3}.\nonumber
             \end{align}

\end{document}